\documentclass{article}
\def\p{\partial}
\def\be{\begin{equation}}
\def\ee{\end{equation}}
\newcommand{\kk}[2]{\frac{#1}{#2}}
\newcommand{\pab}[2]{\kk{\p #1}{\p #2}}
\newcommand{\ff}[1]{{\bf #1}}
\def\a{\alpha}
\def\b{\beta}

\def\e{\epsilon}
\def\lam{\lambda}
\def\s{\qquad}

\def\={\approx}
\def\x{\ff{x}}
\def\vcode#1#2#3{\begin{figure}
\begin{center}
\begin{minipage}[c]{#1\textwidth}
\hrule \vspace{2pt}
{{\bf begin}  \\   %
{\sf #2} \\ {\bf end} \vspace{5pt} \hrule }
\end{minipage}
\caption{#3}
\end{center}   \end{figure}} %%

\begin{document}

\title{Accelerated Particle Swarm Optimization and Support Vector Machine for Business Optimization and Applications }

\author{Xin-She Yang$^1$, Suash Deb$^2$  and Simon Fong$^3$ \\ \\
1) Department of Engineering,
University of Cambridge, \\
Trumpinton Street,
Cambridge CB2 1PZ, UK. \\
\and
2) Department of Computer Science \& Engineering, \\
C. V. Raman College of Engineering, \\
Bidyanagar, Mahura, Janla, 
Bhubaneswar 752054, INDIA. \\
\and
3) Department of Computer and Information Science, \\
Faculty of Science and Technology, \\
University of Macau, Taipa, Macau. \\
}

\date{}

\maketitle

\begin{abstract}

Business optimization is becoming increasingly important because all business activities
aim to maximize the profit and performance of products and services, under
limited resources and appropriate constraints. Recent developments in support vector machine
and metaheuristics show many advantages of these techniques.
In particular, particle swarm optimization is now widely used in solving tough optimization
problems. In this paper, we use a
combination of a recently developed Accelerated PSO and a nonlinear support vector machine to
form a framework for solving business optimization problems.
We first apply the proposed APSO-SVM to production optimization, and then use it for
income prediction and project scheduling.  We also carry out some parametric studies and discuss the
advantages of the proposed metaheuristic SVM.  \\ \\


\noindent Reference to this paper should be made as follows: \\ \\
Yang, X. S., Deb, S., and Fong, S., (2011), Accelerated Particle Swarm Optimization
and Support Vector Machine for Business Optimization and Applications, in: Networked Digital Technologies (NDT2011),
Communications in Computer and Information Science, Vol. 136, Springer, 
pp. 53-66 (2011).

\end{abstract}

%% Begin of the Main Text %%

\section{Introduction}

Many business activities often have to deal with large, complex databases.
This is partly driven by information technology, especially the Internet,
and partly driven by the need to extract meaningful knowledge by data mining.
To extract useful information among a huge amount of data requires
efficient tools for processing vast data sets. This is equivalent to
trying to find an optimal solution to a highly nonlinear problem with multiple, complex
constraints, which is a  challenging task. Various techniques for such data mining and optimization have been
developed over the past few decades. Among these techniques,
support vector machine is one of the best techniques for regression,
classification and data mining \cite{Howley,Kim,Pai,Shi,Shi2,Vapnik}.

On the other hand, metaheuristic algorithms also become powerful for solving
tough nonlinear optimization problems \cite{Blum,Kennedy,Kennedy2,Yang,Yang2}.
Modern metaheuristic algorithms have been developed with an aim to carry out
global search, typical examples are genetic algorithms \cite{Gold},
particle swarm optimisation (PSO) \cite{Kennedy}, and Cuckoo Search \cite{YangDeb,YangDeb2}.
The efficiency of metaheuristic algorithms can be attributed to the
fact that they imitate the best features in nature, especially the selection of the fittest
in biological systems which have evolved by natural selection over millions of years.
Since most data have noise or associated randomness, most these algorithms
cannot be used directly. In this case, some form of averaging or reformulation of the problem
often helps. Even so, most algorithms become difficult to implement for such type of optimization.

In addition to the above challenges, business optimization often concerns
with a large amount but often incomplete data, evolving dynamically over
time. Certain tasks cannot start before other required tasks are completed,
such complex scheduling is often NP-hard and no universally efficient tool exists.
Recent trends indicate that metaheuristics can be very promising, in combination with
other tools such as neural networks and support vector machines \cite{Howley,Kim,Tabu,Smola}.

In this paper, we intend to present a simple framework of business optimization using
a combination of support vector machine with accelerated PSO. The paper is
organized as follows: We first will briefly review particle swarm optimization and accelerated PSO,
and then introduce the basics of support vector machines (SVM).  We then
use three case studies to test the proposed framework. Finally, we discussion its implications
and possible extension for further research.

\section{Accelerated Particle Swarm Optimization}

\subsection{PSO}

Particle swarm optimization (PSO) was developed by Kennedy and
Eberhart in 1995 \cite{Kennedy,Kennedy2}, based on the swarm behaviour such
as fish and bird schooling in nature. Since then, PSO has
generated much wider interests, and forms an exciting, ever-expanding
research subject, called swarm intelligence. PSO has been applied
to almost every area in optimization, computational intelligence,
and design/scheduling applications. There are at least two dozens of
PSO variants, and hybrid algorithms by combining PSO
with other existing algorithms are also increasingly popular.

PSO searches the space of an objective function
by adjusting the trajectories of individual agents,
called particles, as the piecewise paths formed by positional
vectors in a quasi-stochastic manner. The movement of a swarming particle
consists of two major components: a stochastic component and a deterministic component.
Each  particle is attracted toward the position of the current global best
$\ff{g}^*$ and its own best location $\x_i^*$ in history,
while at the same time it has a tendency to move randomly.

Let $\x_i$ and $\ff{v}_i$ be the position vector and velocity for
particle $i$, respectively. The new velocity vector is determined by the
following formula
\be \ff{v}_i^{t+1}= \ff{v}_i^t  + \a \ff{\e}_1
[\ff{g}^*-\x_i^t] + \b \ff{\e}_2  [\x_i^*-\x_i^t].
\label{pso-speed-100}
\ee
where $\ff{\e}_1$ and $\ff{\e}_2$ are two random vectors, and each
entry taking the values between 0 and 1.
The parameters $\a$ and $\b$ are the learning parameters or
acceleration constants, which can typically be taken as, say, $\a\=\b \=2$.

There are many variants which extend the standard PSO
algorithm,  and the most noticeable improvement is probably to use an inertia function $\theta
(t)$ so that $\ff{v}_i^t$ is replaced by $\theta(t) \ff{v}_i^t$
\be \ff{v}_i^{t+1}=\theta \ff{v}_i^t  + \a \ff{\e}_1
[\ff{g}^*-\x_i^t] + \b \ff{\e}_2  [\x_i^*-\x_i^t],
\label{pso-speed-150}
\ee
where $\theta \in (0,1)$ \cite{Chat,Clerc}. In the simplest case,
the inertia function can be taken as a constant, typically $\theta \= 0.5 \sim 0.9$.
This is equivalent to introducing a virtual mass to stabilize the motion
of the particles, and thus the algorithm is expected to converge more quickly.

\subsection{Accelerated PSO}

The standard particle swarm optimization uses both the current global best
$\ff{g}^*$ and the individual best $\ff{x}^*_i$. The reason of using the individual
best is primarily to increase the diversity in the quality solutions, however,
this diversity can be simulated using some randomness. Subsequently, there is
no compelling reason for using the individual best, unless the optimization
problem of interest is highly nonlinear and multimodal.

A simplified version
which could accelerate the convergence of the algorithm is to use the global
best only. Thus, in the accelerated particle swarm optimization (APSO) \cite{Yang,Yang2}, the
velocity vector is generated by a simpler formula
\be \ff{v}_i^{t+1}=\ff{v}_i^t + \a \ff{\e}_n + \b (\ff{g}^*-\x_i^t), \label{pso-sys-10} \ee
where  $\ff{\e}_n$ is drawn from $N(0,1)$
to replace the second term.
The update of the position
is simply \be \x_i^{t+1}=\x_i^t + \ff{v}_i^{t+1}.  \label{pso-sys-20} \ee In order to
increase the convergence even further, we can also write the
update of the location  in a single step
\be  \x_i^{t+1}=(1-\b) \x_i^t+\b \ff{g}^* +\a \ff{\e}_n.  \label{APSO-500} \ee
This simpler version will give the same order of convergence.
Typically, $\alpha = 0.1 L \sim 0.5 L$ where $L$ is the scale of each variable, while $\beta= 0.1 \sim 0.7$
is sufficient for most applications. It is worth pointing out that
 velocity does not appear in equation (\ref{APSO-500}), and there is no need to deal with
 initialization of velocity vectors.
 Therefore, APSO is much simpler. Comparing with many PSO variants, APSO uses only two parameters,
 and the mechanism is simple to understand.

A further improvement to the accelerated PSO is to reduce the randomness
as iterations proceed.
This means   that we can use a monotonically decreasing function such as
\be \a =\a_0 e^{-\gamma t}, \ee
or \be \a=\a_0 \gamma^t, \s (0<\gamma<1), \ee
where $\a_0 \=0.5 \sim 1$ is the initial value of the randomness parameter.
Here $t$ is the number of iterations or time steps.
$0<\gamma<1$ is a control parameter \cite{Yang2}. For example, in our implementation, we will use
\be \a=0.7^t, \ee
where $t \in [0,t_{\max}]$ and $t_{\max}$ is the maximum of iterations.

\section{Support Vector Machine}

Support vector machine (SVM) is an efficient tool for data mining and classification
\cite{Vapnik2,Vapnik3}.  Due to the vast volumes of
data in business, especially e-commerce, efficient
use of data mining techniques becomes a necessity.
In fact, SVM can also be considered as an optimization tool, as its objective is
to maximize the separation margins between data sets. The proper combination of SVM
with metaheuristics could be advantageous.

\subsection{Support Vector Machine}

A support vector machine essentially transforms a set of data
into a significantly higher-dimensional space by nonlinear transformations
so that regression and data fitting can be carried out in this high-dimensional space.
This methodology can be used for data classification, pattern recognition,
and regression, and its theory was based on statistical machine learning theory
\cite{Smola,Vapnik,Vapnik2}.

For classifications with  the learning examples
or data  $(\x_i, y_i)$ where $i=1,2,...,n$ and $y_i \in \{-1,+1\}$,
the aim of the learning is to find a function
$\phi_{\a} (\x)$ from allowable functions $\{\phi_{\a}: \a \in \Omega \}$
such that $\phi_{\a}(\x_i) \mapsto y_i$ for  $(i=1,2,...,n)$
and that the expected risk $E(\a)$ is minimal. That is the minimization
of the risk \be E(\a)=\kk{1}{2} \int |\phi_{\a}(x) -y| dQ(\x,y), \ee
where $Q(\x,y)$ is an unknown probability distribution, which makes
it impossible to calculate $E(\a)$ directly. A simple approach is
to use the so-called empirical risk
\be E_p(\a) \= \kk{1}{2 n} \sum_{i=1}^n  \big|\phi_{\a}(\x_i)-y_i \big|. \ee
However, the main flaw of this approach is that a small risk or error on
the training set does not necessarily guarantee a small
error on prediction if the number $n$ of training data is small \cite{Vapnik3}.

For a given probability of at least $1-p$, the Vapnik bound for the
errors can be written as
\be E(\a) \le R_p(\a) + \Psi \Big(\kk{h}{n}, \kk{\log (p)}{n} \Big), \ee
where
\be \Psi \big(\kk{h}{n}, \kk{\log(p)}{n} \big) =\sqrt{\kk{1}{n} \big[h (\log \kk{2n}{h}+1)
-\log(\kk{p}{4})\big]}.  \ee
Here $h$ is a parameter, often referred to as the Vapnik-Chervonenskis
dimension or simply VC-dimension \cite{Vapnik}, which describes the capacity
for prediction of the function set $\phi_{\a}$.

In essence, a linear support vector machine is to
construct two hyperplanes as far away as possible
and no samples should be between these two planes.
Mathematically, this is equivalent to two equations
\be \ff{w} \cdot \x + \ff{b} = \pm 1, \ee
and a main objective of
constructing these two hyperplanes is to maximize the distance (between the two planes)
\be d=\kk{2}{||\ff{w}||}. \ee
Such maximization of $d$ is equivalent to the minimization of $||w||$ or more conveniently $||w||^2$.
From the optimization point of view, the maximization of margins can be written as
\be \textrm{minimize } \kk{1}{2} ||\ff{w}||^2 = \kk{1}{2} (\ff{w} \cdot \ff{w}). \ee
This essentially becomes an optimization problem
\be \textrm{minimize } \Psi= \kk{1}{2} || \ff{w} ||^2 +\lam \sum_{i=1}^n \eta_i, \ee
\be \textrm{subject to } y_i (\ff{w} \cdot \x_i + \ff{b}) \ge 1-\eta_i, \label{svm-ineq-50} \ee
\be \s \s \s \eta_i \ge 0, \s (i=1,2,..., n), \ee
where $\lam>0$ is a parameter to be chosen appropriately.
Here, the term $\sum_{i=1}^n \eta_i$ is essentially a measure of
the upper bound of the number of misclassifications on the training data.

\subsection{Nonlinear SVM and Kernel Tricks}

The so-called kernel trick is an important technique, transforming data dimensions
while simplifying computation.
By using Lagrange multipliers $\a_i \ge 0$, we can rewrite the above constrained optimization
into an unconstrained version, and we have
\be L=\kk{1}{2} ||\ff{w}||^2 +\lam \sum_{i=1}^n \eta_i - \sum_{i=1}^n \a_i [y_i (\ff{w} \cdot \x_i + \ff{b}) -(1-\eta_i)]. \ee
From this, we can write the Karush-Kuhn-Tucker conditions
\be \pab{L}{\ff{w}}=\ff{w} - \sum_{i=1}^n \a_i y_i \x_i =0, \ee
\be \pab{L}{\ff{b}} = -\sum_{i=1}^n \a_i y_i =0, \ee
\be y_i (\ff{w} \cdot \x_i+\ff{b})-(1-\eta_i) \ge 0, \ee
\be \a_i [y_i (\ff{w} \cdot \x_i + \ff{b}) -(1-\eta_i)]=0, \s (i=1,2,...,n), \label{svm-KKT-150} \ee
\be \a_i \ge 0, \s \eta_i \ge 0, \s (i=1,2,...,n). \ee
From the first KKT condition, we get
\be \ff{w}=\sum_{i=1}^n y_i \a_i \x_i. \ee
It is worth pointing out here that
only the nonzero $\a_i$ contribute to overall solution. This comes from the
KKT condition (\ref{svm-KKT-150}),
which implies that when $\a_i \ne 0$, the inequality (\ref{svm-ineq-50}) must be
satisfied exactly, while $\a_0=0$ means the inequality is automatically
met. In this latter case, $\eta_i=0$. Therefore, only the corresponding training
data $(\x_i, y_i)$ with $\a_i>0$ can contribute to the solution, and thus such
$\x_i$ form the support vectors (hence, the name support vector machine). \index{support vectors}
All the other data with $\a_i=0$ become irrelevant.

It has been shown that the solution for $\a_i$ can be found
by solving the following quadratic programming  \cite{Vapnik,Vapnik3}
\be \textrm{maximize } \sum_{i=1}^n \a_i -\kk{1}{2} \sum_{i,j=1}^n \a_i \a_j y_i y_j (\x_i \cdot \x_j), \ee
subject to
\be \sum_{i=1}^n \a_i y_i=0, \s 0 \le \a_i \le \lam, \s (i=1,2,...,n). \ee
From the coefficients $\a_i$, we can write the final classification or decision
function as \be f(\x) =\textrm{sgn} \big [ \sum_{i=1}^n \a_i y_i (\x \cdot \x_i) + \ff{b} \big ], \ee
where sgn is the classic sign function.

As most problems are nonlinear in business applications, and the above linear SVM cannot be
used. Ideally, we should find some nonlinear transformation $\phi$ so that the
data can be mapped onto a high-dimensional space where the classification
becomes linear. The transformation should be
chosen in a certain way so that their dot product leads to a
kernel-style function $K(\x,\x_i)=\phi(\x) \cdot \phi(\x_i)$.
In fact, we do not need to know such transformations,
we can directly use the kernel functions $K(\x,\x_i)$ to complete this task.
This is the so-called kernel function trick. Now the main task is to chose
a suitable kernel function for a given, specific problem.

For most problems in nonlinear support vector machines, we can
use $K(\x,\x_i)=(\x \cdot x_i)^d$ for polynomial classifiers,
$K(\x,\x_i)=\tanh[k (\x \cdot \x_i) +\Theta)]$ for neural networks,
and by far the most widely used kernel is the Gaussian radial basis function (RBF)
\be K(\x,\x_i)=\exp \Big[-\kk{||\x-\x_i||^2}{(2 \sigma^2)} \Big]
 =\exp \Big[-\gamma ||\x-\x_i||^2 \Big],\ee
for the nonlinear classifiers. This kernel can easily be extended to
any high dimensions. Here $\sigma^2$ is the variance and $\gamma=1/2\sigma^2$ is
a constant. In general, a simple bound of $0 < \gamma \le C$ is used, and
here $C$ is a constant.

Following the similar procedure as discussed earlier for linear SVM,
we can obtain the coefficients $\a_i$ by solving the following optimization
problem \be \textrm{maximize } \sum_{i=1}^n \a_i -\kk{1}{2} \a_i \a_j y_i y_j K(\x_i,\x_j). \ee
It is worth pointing out under Mercer's conditions for the kernel function,
the matrix $\ff{A}=y_i y_j K(\x_i, \x_j)$ is a symmetric positive definite matrix \cite{Vapnik3}, which
implies that the above maximization is a quadratic programming problem, and can thus
be solved efficiently by standard QP techniques \cite{Smola}.

\section{Metaheuristic Support Vector Machine with APSO}

\subsection{Metaheuristics}

There are many metaheuristic algorithms for optimization and most these algorithms
are inspired by nature \cite{Yang}. Metaheuristic algorithms such as genetic
algorithms and simulated annealing are widely used, almost routinely, in many
applications, while relatively new algorithms such as particle swarm optimization \cite{Kennedy},
firefly algorithm and cuckoo search are becoming more and more popular \cite{Yang,Yang2}.
Hybridization of these algorithms with existing algorithms are also emerging.

The advantage of such a combination is to use a balanced tradeoff between
global search which is often slow and a fast local search. Such a balance
is important, as highlighted by the analysis by Blum and Roli \cite{Blum}.
Another advantage of this method is that we can use any algorithms we like
at different stages of the search or even at different stage of iterations.
This makes it easy to combine the advantages of various algorithms so as
to produce better results.

Others have attempted to carry out parameter optimization associated with neural
networks and SVM. For example, Liu et al.  have used SVM optimized by PSO for
tax forecasting \cite{Liu}. Lu et al. proposed a model for finding optimal parameters in SVM by PSO
optimization \cite{Lu}. However, here we intend to propose a generic framework
for combining efficient APSO with SVM, which can be extended to
other algorithms such as firefly algorithm \cite{YangFA,Yang2010}.

\subsection{APSO-SVM}

Support vector machine has a major advantage, that is, it is less likely to
overfit, compared with other methods such as regression and neural networks.
In addition,  efficient quadratic programming can be used for training support vector machines.
However, when there is noise in the data, such algorithms are not quite suitable.
In this case, the learning or training to estimate the parameters in the
SVM becomes difficult or inefficient.

Another issue is that the choice of the values of
kernel parameters $C$ and $\sigma^2$ in the kernel functions;
however, there is no agreed guideline on how
to choose them, though the choice of their values should make
the SVM as efficiently as possible. This itself
is essentially an optimization problem.

Taking this idea further, we first use an educated guess set of
values and use the metaheuristic algorithms such as accelerated PSO
or cuckoo search to find the best kernel parameters
such as $C$ and $\sigma^2$ \cite{Yang,YangDeb}.
Then, we used these parameters to construct the support vector machines
which are then used for solving the problem of interest. During the iterations
and optimization, we can also modify kernel parameters and
evolve the SVM accordingly. This framework can be called a metaheuristic
support vector machine. Schematically, this Accelerated PSO-SVM can be represented as
shown in Fig. 1.

\vcode{0.7}{
Define the objective;  \\
Choose  kernel functions; \\
Initialize various parameters; \\
{\bf while} (criterion) \\
\indent $\quad$ Find optimal kernel parameters by APSO; \\
\indent $\quad$ Construct the support vector machine; \\
\indent $\quad$ Search for the optimal solution by APSO-SVM; \\
\indent $\quad$ Increase the iteration counter; \\
{\bf end} \\
Post-processing the results;
}{Metaheuristic APSO-SVM. }

For the optimization of parameters and business applications discussed below, APSO
is used for both local and global search \cite{Yang,Yang2}.

\section{Business Optimization Benchmarks}

Using the framework discussed earlier, we can easily implement it in
any programming language, though we have implemented using Matlab.
We have validated our implementation using the standard test
functions, which confirms the correctness of the implementation.
Now we apply it to carry out  case studies with known analytical solution
or the known optimal solutions. The Cobb-Douglas production
optimization has an analytical solution which can be used for
comparison, while the second case is a standard benchmark in
resource-constrained project scheduling \cite{Kol}.

\subsection{Production Optimization}

Let us first use the proposed approach to study
the classical Cobb-Douglas production optimization.
For a production of a series of products and the labour costs,
the utility function can be written
\be q=\prod_{j=1}^n u_j^{\a_j} =u_1^{\a_1} u_2^{\a_2} \cdots u_n^{\a_n}, \ee
where all exponents $\a_j$ are non-negative, satisfying
\be \sum_{j=1}^n \a_j =1. \ee
The optimization is the minimization of the utility
\be \textrm{minimize } q \ee
\be \textrm{subject to } \sum_{j=1}^n w_j u_j =K, \ee
where $w_j (j=1,2,...,n)$ are known weights.

This problem can be solved using the Lagrange multiplier method as
an unconstrained problem
\be \psi=\prod_{j=1}^n u_j^{\a_j} + \lam (\sum_{j=1}^n w_j u_j -K), \ee
whose optimality conditions are
\be \pab{\psi}{u_j} = \a_j u_j^{-1} \prod_{j=1}^n u_j^{\a_j} + \lam w_j =0, \quad (j=1,2,...,n), \ee
\be \pab{\psi}{\lam} = \sum_{j=1}^n w_j u_j -K =0. \ee
The solutions are
\be u_1=\kk{K}{w_1 [1+\kk{1}{\a_1} \sum_{j=2}^n \a_j ]}, \;
u_j=\kk{w_1 \a_j}{w_j \a_1} u_1, \ee
where $(j=2,3, ..., n)$.
For example, in a special case of $n=2$, $\a_1=2/3$, $\a_2=1/3$, $w_1=5$,
$w_2=2$ and $K=300$, we have
\[ u_1=\kk{Q}{w_1 (1+\a_2/\a_1)} =40, \;
 u_2=\kk{K \a_2}{w_2 \a_1 (1+\a_2/\a_1)}=50. \]

As most real-world problem has some uncertainty, we can now add
some noise to the above problem. For simplicity, we just  modify
the constraint as
\be \sum_{j=1}^n w_j u_j = K (1+ \beta \epsilon), \ee
where $\epsilon$ is a random number drawn from a Gaussian distribution with a zero mean
and a unity variance,
and $0 \le \beta \ll 1$ is a small positive number.

We now solve this problem as an optimization problem by the proposed APSO-SVM.
In the case of $\beta=0.01$,
the results have been summarized in Table 1
where the values are provided with different problem size $n$ with
different numbers of iterations. We can see that the results converge
at the optimal solution very quickly.

\begin{table}[ht]
\caption{Mean deviations from the optimal solutions.}
\centering
\begin{tabular}{lllll}
\hline \hline
size $n$ & Iterations & deviations \\
\hline
10  & 1000 & 0.014 \\
20 & 5000  & 0.037 \\
50 & 5000   & 0.040 \\
50 & 15000  & 0.009 \\
\hline
\end{tabular}
\end{table}

\section{Income Prediction}

Studies to improve the accuracy of classifications are extensive. For example, Kohavi proposed a
decision-tree hybrid in 1996 \cite{UCI}. Furthermore, an efficient training algorithm for support vector machines was proposed by Platt in 1998 \cite{Platt,Platt2},
and it has some significant impact on machine learning, regression and data mining.

A well-known benchmark for classification and regression is the income prediction using the
data sets from a selected 14 attributes of a household from a sensus form \cite{UCI,Platt}.
We use the same data sets at ftp://ftp.ics.uci.edu/pub/machine-learning-databases/adult
for this case study. There are 32561 samples in the training set with 16281 for testing.
The aim is to predict if an individual's income is above or below 50K ?

Among the 14 attributes, a subset can be selected, and a subset such as age, education level,
occupation, gender and working hours are commonly used.

Using the proposed APSO-SVM and choosing the limit value of $C$ as $1.25$,
the best error of $17.23\%$ is obtained (see Table \ref{table-3}), which is comparable with most accurate predictions
reported in \cite{UCI,Platt}.

\begin{table}[ht]
\caption{Income prediction using APSO-SVM. \label{table-3}}
\centering
\begin{tabular}{l|l|l}
\hline \hline
Train set (size) & Prediction set  & Errors (\%) \\
\hline
512 & 256 & $24.9$ \\
1024 & 256 & $20.4$ \\
16400 & 8200 & $17.23$ \\
\hline
\end{tabular}
\end{table}

\subsection{Project Scheduling}

Scheduling is an important class of discrete optimization with a wider range
of applications in business intelligence. For resource-constrained project scheduling problems,
there exists a standard benchmark
library by Kolisch and Sprecher \cite{Kol,Kol2}. The basic model consists
of $J$ activities/tasks, and some activities cannot start before all its predecessors $h$
are completed. In addition, each activity $j=1,2,...,J$
can be carried out, without interruption,
in one of the $M_j$ modes, and performing any activity $j$ in any chosen
mode $m$ takes $d_{jm}$ periods, which is supported by a set of renewable
resource $R$ and non-renewable resources $N$. The project's makespan or upper
bound is T, and the overall capacity of non-renewable resources is
$K_r^{\nu}$ where $r \in N$. For an activity $j$ scheduled in mode $m$,
it uses $k^{\rho}_{jmr}$ units of renewable resources
and $k^{\nu}_{jmr}$ units of non-renewable resources
in period $t=1,2,..., T$.

For activity $j$, the shortest duration is fit into the time
windows  $[EF_j, LF_j]$ where $EF_j$ is the earliest finish times,
and $LF_j$ is the latest finish times. Mathematically, this
model can be written as \cite{Kol}
\be \textrm{Minimize }\; \Psi (\x) \sum_{m=1}^{M_j} \sum_{t=EF_j}^{LF_j} t \cdot x_{jmt}, \ee
subject to
\[ \sum_{m=1}^{M_h} \sum_{t=EF_j}^{LF_j} t \;\;  x_{hmt} \le \sum_{m=1}^{M_j} \sum_{t=EF_j}^{LF_j} (t-d_{jm}) x_{jmt},
(j=2,..., J), \]
\[ \sum_{j=1}^J \sum_{m=1}^{M_j} k^{\rho}_{jmr} \sum_{q=\max\{t,EF_j\}}^{\min\{t+d_{jm}-1,LF_j\}} x_{jmq} \le K_r^{\rho},
 (r \in R), \]
\be \sum_{j=1}^J \sum_{m=1}^{M_j} k_{jmr}^{\nu} \sum_{t=EF_j}^{LF_j} x_{jmt} \le K^{\nu}_r, (r \in N), \ee
and
\be \sum_{j=1}^{M_j} \sum{t=EF_j}^{LF_j} =1, \s j=1,2,...,J, \ee
where $x_{jmt} \in \{0,1\}$ and $t=1,...,T$.
As $x_{jmt}$ only takes two values $0$ or $1$, this problem
can be considered as a classification problem, and metaheuristic
support vector machine can be applied naturally.

\begin{table}[ht]
\caption{Kernel parameters used in SVM.}
\centering
\begin{tabular}{l|l}
\hline \hline
Number of iterations & SVM kernel parameters \\
\hline
1000 & $C=149.2$, $\sigma^2=67.9$ \\
5000 & $C=127.9$, $\sigma^2=64.0$ \\
\hline
\end{tabular}
\end{table}

Using the online benchmark library \cite{Kol2}, we have solved this type
of problem with $J=30$ activities (the standard test set j30).  The run time
on a modern desktop computer is about 2.2 seconds for $N=1000$ iterations
to 15.4 seconds for $N=5000$ iterations. We have run
the simulations for 50 times so as to obtain meaningful statistics.

The optimal kernel parameters found for the support vector machines
are listed in Table 3, while the deviations from the known  best solution
are given in Table 4 where the results by other methods are also compared.

\begin{table}[ht]
\caption{Mean deviations from the optimal solution (J=30).}
\centering
\begin{tabular}{lllll}
\hline \hline
Algorithm & Authors & $N=1000$ & $5000$ \\
\hline
PSO \cite{Tcho} & Kemmoe et al. (2007)  & 0.26 & 0.21 \\
hybribd GA \cite{Valls} & Valls eta al. (2007) & 0.27 & 0.06 \\
Tabu search \cite{Tabu} & Nonobe \& Ibaraki (2002) & 0.46 & 0.16 \\
Adapting GA \cite{Hart} & Hartmann (2002) & 0.38 & 0.22 \\
{\bf Meta APSO-SVM } & this paper & {\bf 0.19 } & {\bf 0.025} \\
\hline
\end{tabular}
\end{table}

From these tables, we can see that the proposed metaheuristic support vector machine
starts very well, and results are comparable with those by other methods such as hybrid
genetic algorithm. In addition, it converges more quickly, as the number of
iterations increases. With the same amount of function evaluations involved,
much better results are obtained, which implies that APSO is very efficient,
and subsequently the APSO-SVM is also efficient in this context. In addition, this
also suggests that this proposed framework is appropriate for automatically choosing
the right parameters for SVM and solving nonlinear optimization problems.

\section{Conclusions}

Both PSO and support vector machines are now widely used as  optimization techniques
in business intelligence. They can also be used for data mining to extract useful information
efficiently. SVM can also be considered as an optimization
technique in many applications including business optimization. When there is noise in data,
some averaging or reformulation may lead to better performance. In addition, metaheuristic
algorithms can be used to find the optimal kernel parameters for a support vector machine
and also to search for the optimal solutions. We have used three very different case studies to demonstrate
such a metaheuristic SVM framework works.

Automatic parameter tuning  and efficiency improvement will be an important topic for
further research. It can be expected that this framework can be used for other applications.
Furthermore, APSO can also be used to combine with other algorithms such as neutral networks
to produce more efficient algorithms \cite{Liu,Lu}. More studies in this area are highly needed.

%% End of text %%

\end{document}